\documentclass[12pt]{amsart}
\usepackage{amscd,amsmath,amssymb,amsfonts}
\theoremstyle{plain}
\newtheorem{thm}{Theorem}
\newtheorem{lem}[thm]{Lemma}
\newtheorem{cor}[thm]{Corollary}

\theoremstyle{definition}

\numberwithin{thm}{section}
\numberwithin{equation}{section}

\newcommand{\sO}{{\mathcal O}}

\newcommand{\C}{{\mathbb C}}

\newcommand{\F}{{\mathbb F}}

\newcommand{\N}{{\mathbb N}}
\renewcommand{\P}{{\mathbb P}}
\newcommand{\Q}{{\mathbb Q}}

\newcommand{\Z}{{\mathbb Z}}

\begin{document}

\title[Chow group and rational point]{
Varieties over a finite field with trivial Chow group
of 0-cycles have a rational point}

\author{H\'el\`ene Esnault}
\address{Mathematik,
Universit\"at Essen, FB6, Mathematik, 45117 Essen, Germany}
\email{esnault@uni-essen.de}
\date{July 20 print, 2002}

\maketitle
\begin{quote}

\end{quote}

\section{Introduction}
Let $X$ be a smooth projective variety of dimension $d$ over a
field $k$. Let $\overline{k(X)}$ be the algebraic closure of its
function field. If the Chow group of 0-cycles $CH_0(X \times_k
\overline{k(X)})$ is equal to $\Z$, then S. Bloch shows in
\cite{B}, Appendix to Lecture 1, that the diagonal $\Delta \in
CH^d(X\times_k X)\otimes_{\Z} \Q$ decomposes. This means there are
a $N\in \N\setminus \{0\} $, a 0-dimensional subscheme $\xi\subset
X$, a divisor $D\subset X$, a dimension $d$ cycle $\Gamma\subset
X\times  D$ such that
\begin{gather} \label{dec}
N\cdot \Delta \equiv \xi\times X + \Gamma.
\end{gather}
 For sake of completeness, we briefly recall  his argument.
Using the norm, one sees that the kernel of $CH_0(X\times_k
k(X))\to CH_0(X\times_k \overline{k(X)})$ is torsion. Thus  up to
torsion, the cycle ${\rm Spec}\ k(X)$ is equivalent in
$CH_0(X\times_k k(X))$ to
 a $\overline{k}$-rational point of $X$, thus, up to torsion,
 to a $\xi$ as above. On the other hand,
$CH_0(X\times_k k(X))$ is the inductive limit of
$CH_0(X\times_k (X\setminus D))$ as $D$ runs over the divisors of $X$.

This  has various consequences on the shape of de Rham or Hodge
cohomologies in characteristic 0 . Let $(\Delta)\in
H^{2d}(X\times X)$ be the cycle class of $\Delta$.  One applies
 the correspondence $[\Delta]_*=p_{2,*}((\Delta)\cup
p_1^*)=[\xi \times X]_* + [\Gamma]_*$ to, for example,
$H^i_{DR}(X)$. Then $[\xi \times X]_*H^i_{DR}(X)=0$ for $i\ge 1$
as the correspondence factors through $H^i_{DR}(\xi)=0$, while $
[\Gamma]_*H^i_{DR}(X)\subset H^i_{DR}(X)$ dies via the restriction
map $H^i_{DR}(X)\to H^i_{DR}(X\setminus D)$. Using the surjection
$H^i_{DR}(X) \to H^i(X,\sO_X)$ to lift classes, and  the
factorization  $H^i_{DR}(X\setminus D)\to H^i(X, \sO_X)$ coming
from Deligne's Hodge theory (\cite{De}), one concludes that
$H^i(X, \sO_X)=0$ for $i\ge 1$. S. Bloch developed this argument,
and variants of it for \'etale cohomology, to kill the algebraic
part of $H^2$ under the representability assumption of the Chow
group of 0-cycles over $\overline{k(X)}$ (Mumford's theorem).

The purpose of this note it to observe that P. Berthelot's rigid
cohomology \cite{Be} has the required properties to make the above
argument work in this framework. If $k$ has characteristic $p >0$,
let $W(k)$ be its ring of Witt vectors, $K$ be the quotient field
of $W(k)$. Let $X$ be smooth proper over $k$, $Z\subset X$ be a
closed subvariety, and $U=(X\setminus Z)$ be its complement. The
rigid cohomology, which coincides with the crystalline cohomology
on $X$,  fulfills a localization sequence (\cite{Be}, 2.3.1)
\begin{gather}\label{loc}
\ldots \to H^i_Z(X/K)\to H^i(X/K)\to H^i(U/K) \to \ldots
\end{gather}
which is compatible with the Frobenius action (\cite{Ch}, Theorem
2.4). By \cite{Bcrys} and \cite{I}, the slope [0\ \ 1[ part of
$H^i(X/K)$ is $H^i(X, W\sO_X)\otimes _{W(k)} K$.
One has
\begin{thm} \label{thm:chow}
Let $X$ be a smooth projective  variety over a perfect  field $k$
of characteristic $p>0$. If the Chow group of 0-cycles $CH_0(X
\times_k \overline{k(X)})$ is equal to $\Z$, then the slope [0 \ \
1[ part of $H^i(X/K)$ is vanishing for $i>0$.
\end{thm}
On the other hand, if one now assumes that $k=\F_q$ is a finite
field, with $q=p^n$, the Lefschetz trace formula for crystalline
cohomology (see e.g. \cite{Et}, section II.1)
\begin{gather} \label{cong}
|X(k)| = \sum_i (-1)^i {\rm Trace} ({\rm Frob}^n|H^i(X/K))
\end{gather}
implies in particular that if all the slopes of Frob are $\ge 1$ and $X$ is
geometrically connected,
then
\begin{gather}
|X(k)|\equiv 1 \ \text{modulo} \ q.
\end{gather}
Thus one has
\begin{cor} \label{cor:point}
Let $X$ be a  smooth, projective, geometrically connected variety
over a finite field $k$. If the Chow group of 0-cycles $CH_0(X
\times_k \overline{k(X)})$ is equal to $\Z$, then $X$ has a
rational point over $k$.
\end{cor}
An example of application is provided by Fano varieties. A variety
$X$ is said to be Fano if it is smooth, projective, geometrically
connected and the dual  of the dualizing sheaf $\omega_X$ is
ample. By \cite{Ko}, Theorem V. 2.13,  Fano varieties on any
algebraically closed field are chain rationally connected
(\cite{Ko}, Definition IV. 3.2) in the sense that any two rational
points can be joined by a chain of rational curves. This implies
in particular that $CH_0(X\times_k \overline{k(X)})=\Z$. Thus one
concludes
\begin{cor} \label{cor:fano}
Let $X$ be a Fano variety over a finite field $k$, or more generally, let
$X$ be a smooth, projective, geometrically connected
variety over a finite field $k$, which is chain rationally connected
over $\overline{k(X)}$.
Then $X$ has a
rational point.
\end{cor}
This corollary answers positively a conjecture by S. Lang
\cite{L} and Yu. Manin \cite{Ma}. This is the reason why we write
down the argument for Theorem \ref{thm:chow}, while it is a direct
adaption of S. Bloch's argument \cite{B} to  crystalline and
rigid cohomologies.

That a study of crystalline cohomology  should yield via the congruence
\eqref{cong} the existence
of a rational point on Fano varieties over finite fields is
entirely due to M. Kim. His idea was to kill the whole cohomology
$H^i(X, W\omega_X)$ for $i< {\rm dim}(X)$, using solely the
structure of crystalline cohomology on $X$, together with its
Verschiebung and Frobenius operators, and using the ampleness of
$\omega_X^{-1}$. The point of Corollary
\ref{cor:point} is that Koll\'ar-Miyaoka-Mori's  and Campana's
theorem (\cite{Ko}, loc. cit.), which is
anchored in geometry, together with Bloch's type Chow group
argument,  force (weaker) cohomological consequences. In a way, the difficult
theorem is the geometric one.

\noindent {\it Acknowledgements}. This note relies on P.
Berthelot's work on rigid cohomology, with which the author is
not familiar. It is a pleasure to thank P. Berthelot,  S. Bloch
and O. Gabber for their substantial help and for their encouragement.
I thank the IHES for support during the preparation of this work.

\section{Proof of theorem \ref{thm:chow}}
In this section we prove Theorem \ref{thm:chow}. Thus $X$ is a
smooth projective variety over $k$. For any codimension $d$ cycle
$Z \subset X\times X$, the correspondence
\begin{gather}
[Z]_*=p_{2,*}((Z)\cup p_1^*)
\end{gather}
is well defined on $H^i(X/K)$. One needs for this the existence
of the cycle class
\begin{gather}
(Z) \in H^{2d}((X\times X)/K)
\end{gather}
which is provided by \cite{Be}, Corollaire 5.7, (ii), and by \cite{P}, 6.2,
for the factorization through the Chow group, the
contravariance for $p_1^*$ and the covariance for $p_{2,*}$,
applied to crystalline cohomology of smooth proper varieties. In
particular, this correspondence  factors through $H^i(\xi/K)$ if
$Z=(\xi \times X)$, which shows via formula \eqref{dec} that
\begin{gather}
N[\Delta]_*= [\Gamma]_*
\end{gather}
on $H^i(X/K)$ for $i>0$. On the other hand, since $\Gamma \subset X\times
D$, $[\Gamma]_*$, seen as a correspondence of $X$ to $(X\setminus
D)$, is trivial.
One has
\begin{gather}
[\Gamma]_*(H^i(X/K))\subset {\rm Ker}(H^i(X/K)\to H^i((X\setminus
D)/K)\\
= \big({\rm by} \ (\ref{loc}) \big)\ {\rm Im}
(H^i_D(X/K)\subset H^i(X/K)) .\notag
\end{gather}
Thus, using \cite{Ch}, Theorem 2.4,   Theorem \ref{thm:chow} is a
consequence of the following
\begin{lem}(P. Berthelot) Let $X$ be a smooth, geometrically connected,
quasi-compact and separated
scheme over a perfect field $k$ of characteristic $p>0$  and let $Z\subset X$ be a non-empty
subvariety of codimension $\ge 1$. Then the slopes of
$H^i_Z(X/K)$ are $\ge 1$.
\end{lem}
\begin{proof}
Let $\ldots \subset Z_i\subset Z_{i-1} \subset \ldots \subset Z_0=
Z$ be a finite stratification by closed subsets such that
$Z_{i-1}\setminus Z_i$ is smooth.   The  localization \cite{Be},
2.5.1
\begin{gather}
\ldots \to H^i_{Z_i} (X/K)\to H^i_{Z_{i-1}}(X/K)\to
H^i_{(Z_{i-1}\setminus Z_i)} ((X\setminus Z_i)/K)\to \ldots
\end{gather}
allows to reduce to the case where $Z$ is smooth. If $X$ is
affine, then the Gysin isomorphism (purity) $H^{i-2\cdot {\rm
codim}(Z)}(Z) \xrightarrow{\cong} H^i_Z(X)$ commutes to $p^{{\rm
codim}(Z)}\cdot$ Frob on $H^{i-2\cdot{\rm codim}(Z)}(Z/K)$ and
Frob on $H^i_Z(X/K)$ (\cite{Ch}, Theorem 2.4). Since the slopes
on $H^i(Z/K)$ are all $\ge 0$, we conclude that the slopes of the
cohomology with support are $\ge 1$. In general, one considers a
finite affine covering $X=\cup_{i=0}^N U_i$. The  spectral
sequence
\begin{gather}
E_2^{ab}=H^a(\ldots \to H^b_Z(U^{a-1}/K)\to H^b_Z(U^a/K) \to H^b_Z(U^{a+1}/K) \to \ldots)
\end{gather}
converges to $H^{a+b}(X/K)$.
Here the open sets $U^a$ are the $(a+1)$ by $(a+1)$ intersections of the
$U_i$ and the maps are the restriction maps. If $H^b_Z(U^a/K)\neq 0$, then $Z$ meets $U^a$ and its slopes are $\ge 1$ by the previous case.
Thus the spectral sequence has only contributions with $\ge 1$ slopes.
 This finishes the proof.
\end{proof}

\section{Comments}
Corollary \ref{cor:fano} can be compared to the main theorems of
\cite{Ha}, resp. \cite{dJ}, where the finite field is replaced by
$k=F({\rm curve})$ for $F=\C$, resp. an algebraically closed field
in characteristic $p>0$. There the authors show the existence of
a rational point on a smooth Fano variety defined over $k$. In
the latter case, one has to add ``separably'' to ``chain
rationally connected''.

On the other hand, if $X$ is a hypersurface $\subset \P^n$  of
degree $\le n$ over a field $k$ of characteristic 0, then $CH_0(X
\times_k \overline{k(X)})=\Z$ by Roitman's theorem (\cite{Ro}),
whether $X$ is smooth or not. It suggests that perhaps there is a
stronger version of Theorem \ref{thm:chow}, requiring $X$
projective but not smooth. This would be compatible with the
congruence results of Ax  and Katz \cite{Ka}, and their Hodge
theoritic counter-part (\cite{DD}, \cite{E}, \cite{ENS}). In this
case it
 says  $H^i(X',
\sO_{X'})=0$, where $X'=\pi^{-1}(X)_{{\rm red}}$ and $\pi: \P\to
\P^n$ is a birational map with $\P$ smooth, such that $X'$ is a
normal crossings divisor. The method used here does not apply to
the singular situation.

After receiving this note, G. Faltings and O. Gabber explained to
me that one could use a similar argument in \'etale cohomology in
order to obtain the same conclusion, and M. Kim replaced the use
of rigid cohomology in the localisation argument presented here by
the one of de Rham-Witt cohomology with logarithmic poles. Since
it does not yield a stronger result, we do not develop their
arguments.

Finally let us observe that, replacing $D$ in the argument by a
subvariety of codimension $\kappa\ge 1$, replaces the congruence
$(1 \ \text{mod} \ q)$ by $(1 \ \text{mod} \ q^\kappa)$. But it is
hard to understand what are the geometric conditions which
guarantee higher codimension. Again, hypersurfaces of very low
degree fulfill the correct congruence (\cite{Ka}), but I do not
know whether it is reflected by this strong Chow group property.
Without this,  the method presented here gives a different proof
of the main result of  \cite{Ka} in the smooth case when
$\kappa=1$.

\bibliographystyle{plain}

\renewcommand\refname{References}

\end{document}